\documentclass[number,citesort,seceqn,dvips]{arxbj}
\usepackage{graphicx}

\aid{0}
\volume{17}
\issue{4}
\pubyear{2011}
\firstpage{1327}
\lastpage{1343}
\doi{10.3150/10-BEJ314}

\makeatletter

\newtheorem{theorem}{Theorem}[section]
\newtheorem{lemma}[theorem]{Lemma}
\newtheorem{cor}[theorem]{Corollary}

\def\eqref#1{(\ref{#1})}

\makeatother

\begin{document}
\begin{frontmatter}

\title{Sharp maximal inequalities for the moments of martingales and
non-negative submartingales}
\runtitle{Maximal inequalities}

\begin{aug}
\author{\fnms{Adam} \snm{Os\textup{\c{E}}kowski}\ead[label=e1]{ados@mimuw.edu.pl}} 
\runauthor{A. Os\c{e}kowski}
\address{Department of Mathematics, Informatics and Mechanics,
University of Warsaw,
Banacha 2, 02-097 Warsaw,
Poland. \printead{e1}}
\pdfauthor{Adam Osekowski}
\end{aug}

\received{\smonth{2} \syear{2010}}

%
\begin{abstract}
In the paper we study sharp maximal inequalities for martingales and
non-negative submartingales: if $f$, $g$ are martingales satisfying
\[
|\mathrm{d}g_n|\leq|\mathrm{d}f_n|,\qquad n=0, 1, 2, \ldots,
\]
almost surely, then
\[
\Bigl\|\sup_{n\geq0} |g_n|\Bigr\|_p \leq p \|f\|_p,\qquad p\geq2,
\]
and the inequality is sharp. Furthermore, if $\alpha\in[0,1]$, $f$ is
a non-negative submartingale and $g$ satisfies
\[
|\mathrm{d}g_n|\leq|\mathrm{d}f_n|\quad \mbox{and}\quad |\mathbb{E}(\mathrm{d}g_{n+1}|\mathcal
{F}_n)|\leq\alpha\mathbb{E}
(\mathrm{d}f_{n+1}|\mathcal{F}_n),\qquad n=0, 1, 2, \ldots,
\]
almost surely, then
\[
\Bigl\|\sup_{n\geq0} |g_n|\Bigr\|_p \leq(\alpha+1)p \|f\|_p,\qquad p\geq2,
\]
and the inequality is sharp. As an application, we establish related
estimates for stochastic integrals and It\^{o} processes.
The inequalities strengthen the earlier classical results of Burkholder
and Choi.

\end{abstract}

%

\end{frontmatter}

\section{Introduction}

The purpose of the paper is to provide the best constants in some
maximal inequalities for martingales and non-negative submartingales.
Let us start with introducing the necessary notation. Let $(\Omega
,\mathcal{F}
,\mathbb{P})$ be a non-atomic probability space, equipped with a
filtration $(\mathcal{F}_n)_{n\geq0}$, that is, a~non-decreasing
family of
sub-$\sigma$-fields of $\mathcal{F}$. Let $f=(f_n)$ and $g=(g_n)$ be adapted,
real-valued integrable processes. The difference sequences $\mathrm{d}f=(\mathrm{d}f_n)$
and $\mathrm{d}g=(\mathrm{d}g_n)$ of $f$ and $g$ are defined by the equations
\[
f_n=\sum_{k=0}^n \mathrm{d}f_k,\qquad g_n=\sum_{k=0}^n \mathrm{d}g_k,\qquad n=0, 1, 2,
\ldots.
\]
We are particularly interested in those pairs $(f,g)$ for which a
certain domination relation is satisfied. Following Burkholder \cite
{B0}, we say that $g$ is \emph{differentially subordinate} to $f$ if,
for any $n\geq0,$ we have
\[
\mathbb{P}(|\mathrm{d}g_n|\leq|\mathrm{d}f_n|)=1.
\]
As an example, let $g$ be a transform of $f$ by a predictable sequence
$v=(v_n)$ bounded in absolute value by $1$; that is, we have $\mathbb
{P}(|v_n|\leq1)=1$ and $\mathrm{d}f_n=v_n\,\mathrm{d}g_n$, $n\geq0$. Here, by
predictability, we mean that $v_0$ is $\mathcal{F}_0$-measurable and
$v_n$ is
$\mathcal{F}_{n-1}$-measurable for $n\geq1$. In the particular case
when each
$v_n$ is deterministic and takes values in $\{-1,1\}$, we will say that
$g$ is a $\pm1$ transform of $f$.

Another domination we will consider is the so-called $\alpha$-strong
subordination, where $\alpha$ is a fixed non-negative number. This
notion was introduced by Burkholder in \cite{B2} in the special case
$\alpha=1$ and extended to a general case by Choi \cite{C}: The
process $g$ is $\alpha$-strongly subordinate to $f$ if it is
differentially subordinate to $f$ and, for any $n\geq0$,
\[
|\mathbb{E}(\mathrm{d}g_{n+1}|\mathcal{F}_{n})|\leq\alpha|\mathbb
{E}(\mathrm{d}f_{n+1}|\mathcal{F}_{n})|
\]
almost surely.

There is a vast literature concerning the comparison of the sizes of
$f$ and $g$ under the assumption of one of the dominations above and
the further condition that $f$ is a martingale or non-negative
submartingale; we refer the interested reader to the papers \cite{B0,B05,B2,C,H,O0,O11,O3,S,W}
and the references therein. In addition, these
inequalities have found their applications in many areas of
mathematics: Banach space theory \cite{Bo,B30}; harmonic analysis
\cite{B15,C2,C3}; functional analysis \cite{B0,B044,S};
analysis \cite{BB,BW}; stochastic integration \cite{B0,B1,O,S,W};
and more. To present our motivation, we state here only two theorems.
Let us start with a fundamental result of Burkholder \cite{B0}. We use
the notation $\|f\|_p=\sup_n \|f_n\|_p$, $p\in[1,\infty]$.

\begin{theorem}[(Burkholder)]\label{burkhr}
Assume that $f$, $g$ are martingales and $g$ is differentially
subordinate to $f$. Then, for any $1<p<\infty$,
%
\begin{equation}\label{strongtype}
\|g\|_p\leq(p^*-1)\|f\|_p,
\end{equation}
where $p^*=\max\{p,p/(p-1)\}$. The constant $p^*-1$ is the best
possible; it is already the best possible if $g$ is assumed to be a
$\pm1$ transform of $f$.
\end{theorem}

Here, by the optimality of the constant, we mean that for any $r<p^*-1$
there exists a martingale $f$ and its $\pm1$ transform $g$, for which
$\|g\|_p>r\|f\|_p$.

The submartingale version of the estimate above is the following result
of Choi \cite{C}.

\begin{theorem}[(Choi)]\label{choi}
Assume that $f$ is a non-negative submartingale and $g$ is $\alpha
$-differentially subordinate to $f$, $\alpha\in[0,1]$. Then for any
$1<p<\infty$,
%
\begin{equation}\label{strongtype2}
\|g\|_p\leq(p_\alpha^*-1)\|f\|_p,
\end{equation}
where $p_\alpha^*=\max\{(\alpha+1)p,p/(p-1)\}$. The constant is the
best possible.
\end{theorem}

In the paper we deal with a considerably harder problem and determine
the optimal constants in the related moment estimates involving the
\emph{maximal functions} of $f$ and $g$. For $n\geq0$, let
$f_n^*=\sup_{0\leq k \leq n}|f_k|$ and $f^*=\sup_{k\geq0} |f_k|$.
Here is our first main result.

\begin{theorem}
Let $f$, $g$ be martingales with $g$ being differentially subordinate
to $f$. Then for any $p\geq2$,
%
\begin{equation}\label{mainin0}
\|g^*\|_p \leq p\|f\|_p
\end{equation}
and the constant $p$ is the best possible. It is already the best
possible in the following weaker inequality: If $f$ is a martingale and
$g$ is its $\pm1$ transform, then
%
\begin{equation}\label{mainin000}
\|g^* \|_p \leq p \| f^*\|_p.
\end{equation}
\end{theorem}

Note that the validity of the estimates \eqref{mainin0} and \eqref
{mainin000} is an immediate consequence of \eqref{strongtype} and
Doob's bound $\|f^*\|_p \leq\frac{p}{p-1}\|f\|_p$, $p> 1$. The
non-trivial (and quite surprising) part is the optimality of the
constant $p$.

Now let us state the submartingale version of the theorem above.

\begin{theorem} \label{main}
Fix $\alpha\in[0,1]$. Let $f$ be a non-negative submartingale and $g$
be real valued and $\alpha$-strongly subordinate to $f$. Then for any
$p\geq2$,
%
\begin{equation}\label{mainin}
\|g^*\|_p \leq(\alpha+1)p\|f\|_p
\end{equation}
and the constant $(\alpha+1)p$ is the best possible. It is already the
best possible in the weaker estimate
%
\begin{equation}\label{mainin0000}
\|g^*\|_p \leq(\alpha+1)p\|f^*\|_p.
\end{equation}
\end{theorem}

There is a natural question: What is the best constant in the
inequalities above in the case $1<p<2$? Unfortunately, we have been
unable to answer it; our reasoning works only for the case $p\geq2$.

The proof of \eqref{mainin} is based on a technique invented by
Burkholder in \cite{B1}. It enables us to translate the problem of
proving a maximal inequality for martingales to that of finding a
certain special function, an upper solution to a corresponding
nonlinear problem. The method can be easily extended to the
submartingale setting (see \cite{O}) and we construct the function in
Section~\ref{sec3}. For the sake of construction, we need a solution to a
differential equation that is analyzed in Section~\ref{sec2}. The next two
sections are devoted to the proofs of the announced results:
Section~\ref{sec4}
contains the proof of the estimate \eqref{mainin} and the final part
concerns the optimality of the constants appearing in \eqref
{mainin000} and \eqref{mainin0000}. In the final section, we present
some applications: sharp estimates for stochastic integrals and It\^{o}
processes.


\section{A differential equation}\label{sec2}

For a fixed $\alpha\in(0,1]$ and $p\geq2$, let $C=C_{p,\alpha
}=[(\alpha+1)p]^p(p-1)$. A central role in the paper is played by a
certain solution to the differential equation
%
\begin{equation}\label{diff}
\gamma'(x)=\frac{-1+C(1-\gamma(x))\gamma(x)x^{p-2}}{1+C(1-\gamma
(x))x^{p-1}}.
\end{equation}

\begin{lemma}
There is a solution $\gamma\dvtx [((\alpha+1)p)^{-1}, \infty) \to\mathbb
{R}$ of
\eqref{diff},
satisfying the initial condition
%
\begin{equation}\label{initial}
\gamma\biggl(\frac{1}{(\alpha+1)p}\biggr)=1-[(\alpha+1)p]^{-1}.
\end{equation}
The solution is non-decreasing, concave and bounded from above by $1$.
\end{lemma}

\begin{pf}
Let $\gamma$ be a solution to \eqref{diff}, satisfying \eqref
{initial} and extended to a maximal subinterval $I$ of $[((\alpha
+1)p)^{-1},\infty)$. It is convenient to split the proof into a few steps.

\emph{Step} 1: $I=[((\alpha+1)p)^{-1},\infty)$. In view of the
Picard--Lindel\"{o}f theorem, this will be established if we show that
$\gamma<1$ on $I$. To this end, suppose that the set $\{x\in I\dvtx \gamma
(x)=1\}$ is non-empty and let $y$ denote its smallest element. Then, by
\eqref{diff}, we have $\gamma'(y)=-1$, which, by minimality of $y$,
implies $\gamma(((\alpha+1)p)^{-1})>1$ and contradicts \eqref{initial}.

\emph{Step \textup{2:} Concavity of $\gamma$.} Suppose that the set $\{x\in I\dvtx
\gamma''(x)>0\}$ is non-empty and let $z$ denote its infimum. Consider
the functions $F,G\dvtx (((\alpha+1)p)^{-1},\infty) \to\mathbb{R}$ given by
\begin{eqnarray*}
F(x)&=&\gamma(x)-x\gamma'(x),
\\
G(x)&=&\bigl(1-\gamma(x)\bigr)x^{p-2}.
\end{eqnarray*}
Observe that
%
\begin{equation}\label{echu3}
G>0 \qquad\mbox{on } I \quad\mbox{and}\quad F>0
\qquad\mbox{on } \bigl(\bigl((\alpha+1)p\bigr)^{-1},z+\varepsilon\bigr)
\end{equation}
for some $\varepsilon>0$. The statement about $G$ is clear, while the
positivity of $F$ follows from
\[
F'(x)=-x\gamma''(x)\geq0,\qquad x\in\bigl(\bigl((\alpha+1)p\bigr)^{-1},z\bigr]
\]
and
\[
F\bigl(\bigl((\alpha+1)p\bigr)^{-1}+\bigr)=\frac{1}{p}>0.
\]
Now multiply \eqref{diff} throughout by $1+C(1-\gamma(x))x^{p-1}$ and
differentiate both sides. We obtain an equality that is equivalent to
%
\begin{equation}\label{echu}
\gamma''(x)\bigl(1+CxG(x)\bigr)=CF(x)G'(x),\qquad x>\frac{1}{(\alpha+1)p}.
\end{equation}
As a first consequence, we have $z>((\alpha+1)p)^{-1}$. To see this,
tend with $x$ down to $((\alpha+1)p)^{-1}$ and observe that $F$ and
$G$ have strictly positive limits; furthermore,
%
\begin{equation}\label{echu4}
G'(x)=x^{p-3}\bigl[(p-2)\bigl(1-\gamma(x)\bigr)-x\gamma'(x)\bigr]=:x^{p-3}J(x)
\end{equation}
with $J(((\alpha+1)p)^{-1})=-\frac{\alpha(p-1)}{(\alpha+1)p}<0$.
Combining \eqref{echu3} and \eqref{echu} we see that, for some
$\varepsilon>0$, $G'\leq0$ on $(z-\varepsilon,z)$ and $G'>0$ on
$(z,z+\varepsilon)$. Consequently, by \eqref{echu4}, $J\leq0$ on
$(z-\varepsilon,z)$ and $J>0$ on $(z,z+\varepsilon)$. This implies
$J'(z)>0$ and since $J'(z)=-(p-1)\gamma'(z)$, we get $\gamma'(z)<0$.
However, this contradicts $G'(z)=0$, in view of \eqref{echu4} and
$\gamma(z)<1$. Let us stress that here, in the last passage, we use
the inequality $p\geq2$.

\emph{Step \textup{3:} $\gamma$ is non-decreasing.} It follows from \eqref
{echu}, the concavity of $\gamma$ and positivity of $F$ and $G$, that
$G'\leq0$, or, by \eqref{echu4},
%
\begin{equation}\label{wazne}
(p-2)\bigl(1-\gamma(x)\bigr)-x\gamma'(x)\leq0.
\end{equation}
The claim follows.
\end{pf}

\begin{figure}

\includegraphics{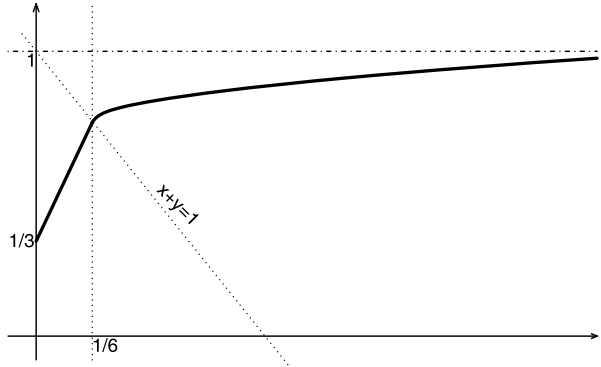}

\caption{The graph of $\gamma$ (the bold line) in the case $p=3$, $\alpha=1$.
Note that $\gamma$ is linear on $[0,1/6]$ and solves \protect\eqref{diff} on
$(1/6,\infty)$.}\label{fig1}
\end{figure}

Let us extend $\gamma$ to the whole half-line $[0,\infty)$ by
\[
\gamma(x)=[(p-1)(\alpha+1)-1]x+\frac{1}{p},\qquad x\in\biggl[0,\frac
{1}{(\alpha+1)p}\biggr).
\]
It can be verified readily that $\gamma$ is of class $C^1$ on
$(0,\infty)$. For the sake of reader's convenience, the graph of $\gamma$, corresponding to $p = 3$
and $\alpha = 1$, is presented on Figure \ref{fig1}.

Let $H\dvtx [((\alpha+1)p)^{-1},\infty) \to[1,\infty)$ be given by
$H(x)=x+\gamma(x)$ and let $h$ be the inverse to~$H$. Clearly, we have
%
\begin{equation}\label{bouh}
x-1 \leq h(x) \leq x,\qquad x\geq1.
\end{equation}
We conclude this section by providing a formula for $h'$ to be used
later. As
%
\begin{equation}\label{backs}
h'(x)=\frac{1}{H'(h(x))}=\frac{1}{1+\gamma'(h(x))},\qquad x>1,
\end{equation}
it can be derived that, in view of \eqref{diff},
%
\begin{equation}\label{prod}
h'(x)=\frac{1+((\alpha+1)p)^p(p-1)(h(x)-x+1)h(x)^{p-1}}{((\alpha
+1)p)^p(p-1)(h(x)-x+1)h(x)^{p-2}x}.
\end{equation}

\section{The special function}\label{sec3}

Throughout this section, $\alpha\in(0,1]$ and $p \geq2$ are fixed.
Let $S$ denote the strip $[0,\infty)\times[-1,1]$. Consider the
following subsets of $S$.
\begin{eqnarray*}
D_0&=&\{(x,y)\in S\dvtx |y|\leq\gamma(x)\},
\\
D_1&=&\{(x,y)\in S\dvtx |y|>\gamma(x), x+|y|\leq1\},
\\
D_2&=&\{(x,y)\in S\dvtx |y|>\gamma(x), x+|y|> 1\}.
\end{eqnarray*}

Introduce the function $u\dvtx S \to\mathbb{R}$ by
\[
u(x,y)=
\cases{
1-[(\alpha+1)p]^px^p &\quad \mbox{on }$ D_0$,\cr
1-\biggl(\dfrac{ px+p|y|-1}{p-1}\biggr)^{p-1}\bigl[p\bigl(p(\alpha+1)-1\bigr)x-p|y|+1\bigr] &\quad \mbox
{on }$ D_1$,\cr
1-[(\alpha+1)p]^ph(x+|y|)^{p-1}[px-(p-1)h(x+|y|)] &\quad \mbox{on }$ D_2$.
}
\]
Let $U\dvtx[0,\infty)\times\mathbb{R}\times(0,\infty) \to\mathbb{R}$
be given by
\[
U(x,y,z)=(|y|\vee z)^pu\biggl(\frac{x}{|y|\vee z},\frac{y}{|y|\vee z}\biggr).
\]

As we will see below, the function $U$ is the key to the inequality
\eqref{mainin}. Let us study the properties of this function.

\begin{lemma}\label{boundy}
The function $U$ is of class $C^1$. Furthermore, there exists an
absolute constant $K$ such that, for all $x> 0, y\in\mathbb{R}, z>0,$
we have
%
\begin{equation}\label{1bou}
U(x,y,z)\leq K(x+|y|+z)^p
\end{equation}
and
%
\begin{equation}\label{2bou}
U_x(x,y,z)\leq K(x+|y|+z)^{p-1}, \qquad U_x(x,y,z)\leq K(x+|y|+z)^{p-1}.
\end{equation}
\end{lemma}

\begin{pf}
The continuity of the partial derivatives can be verified readily. The
inequality \eqref{1bou} is evident for those $(x,y,z)$, for which
$(\frac{x}{|y|\vee z},\frac{y}{|y|\vee z})\in D_0\cup D_1$; for the
remaining $(x,y,z),$ it suffices to use \eqref{bouh}. Finally, the
inequality \eqref{2bou} is clear if $(\frac{x}{|y|\vee z},\frac
{y}{|y|\vee z})\in D_0\cup D_1$. For the remaining points one applies
\eqref{bouh} and \eqref{backs}, the latter inequality implying $h'<1$.
\end{pf}

Now let us deal with the following majorization property.

\begin{lemma}
For any $(x,y,z)\in[0,\infty)\times\mathbb{R}\times(0,\infty)$,
we have
%
\begin{equation}\label{maj}
U(x,y,z) \geq(|y| \vee z)^p-[(\alpha+1)p]^px^p.
\end{equation}
\end{lemma}

\begin{pf}
The inequality is equivalent to $u(x,y)\geq1-[(\alpha+1)p]^px^p$ and
we need to establish it only on $D_1$ and $D_2$. On $D_1$, the
substitutions $X=px$ and $Y=p|y|-1$ (note that $Y\geq0$) transform it into
\[
(\alpha+1)^pX^p\geq\biggl({ \frac{ X+Y}{p-1}}\biggr)^{p-1}\bigl[\bigl(p(\alpha+1)-1\bigr)X-Y\bigr].
\]
This inequality is valid for all non-negative $X$, $Y$. To see this,
observe that by homogeneity we may assume $X+Y=1$, and then the
estimate reads
\[
F(X):=(\alpha+1)^pX^p- (p-1)^{-p+1}[p(\alpha+1)X-1]\geq0,\qquad X\in[0,1].
\]
Now it suffices to note that $F$ is convex on $[0,1]$ and satisfies
\[
F\biggl(\frac{1}{(p-1)(\alpha+1)}\biggr)=F'\biggl(\frac{1}{(p-1)(\alpha+1)}\biggr)=0.
\]

It remains to show the majorization on $D_2$. It is dealt with in a
similar manner: Setting $s=x+|y|>1$, we see that \eqref{maj} is
equivalent to
\[
G(x):=x^p-h(s)^{p-1}[px-(p-1)h(s)]\geq0,\qquad s-1<x<h(s).
\]
It is easily verified that $G$ is convex and satisfies
$G(h(s))=G'(h(s))=0$. This completes the proof of \eqref{maj}.
\end{pf}

The main property of the function $U$ is the concavity along the lines
of slope belonging to $[-1,1]$.

\begin{lemma}\label{conv}
For fixed $y, z$ satisfying $z>0$, $|y|\leq z$, and any $a\in[-1,1]$,
the function $\Phi=\Phi_{y,z,a}\dvtx [0,\infty) \to\mathbb{R}$ given by
\[
\Phi(t)=U(t,y+at,z)
\]
is concave.
\end{lemma}

Before we turn to the proof, let us first establish some useful consequences.

\begin{cor} \label{coro3}
\textup{(i)} The function $U$ has the following property: For any $x, y, z, k_x,
k_y$ such that $x, x+k_x\geq0$, $z>0$, $|y|\leq z$ and $|k_y|\leq
|k_x|,$ we have
%
\begin{equation}\label{convo}
U(x+k_x,y+k_y,z) \leq U(x,y,z)+U_x(x,y,z)k_x+U_y(x,y,z)k_y
\end{equation}
(for $x=0,$ we replace $U_x(0,y,z)$ by right-sided derivative $U_x(0+,y,z)$).

\textup{(ii)} For any $x\geq1,$ we have
%
\begin{equation}\label{start}
U(x,1,1)\leq0.
\end{equation}
\end{cor}

\begin{pf} (i) This follows immediately.

(ii) We have $\Phi_{0,1,x^{-1}}(0)=U(0,0,1)=1$ and $\Phi
_{0,1,x^{-1}}(((\alpha+1)p)^{-1})=U(((\alpha+1)p)^{-1},\break x^{-1}((\alpha
+1)p)^{-1},1)=0$, since $(((\alpha+1)p)^{-1},x^{-1}((\alpha
+1)p)^{-1},1) \in D_0$. Since $x\geq1>((\alpha+1)p)^{-1}$, the lemma
above gives
$U(x,1,1)=\Phi_{0,1,x^{-1}}(x)\leq0.$
\end{pf}

\begin{pf*}{Proof of Lemma \ref{conv}}
By homogeneity, we may assume $z=1$. As $\Phi$ is of class $C^1$, it
suffices to verify that $\Phi''(t)\leq0$ for those $t$, for which
$(t,y+at)$ lies in the interior of $D_0$, $D_1$, $D_2$ or outside the
strip $S$. Since $U(x,y,z)=U(x,-y,z)$, we may restrict ourselves to the
case $y+at\geq0$. If $(t,y+at)$ belongs to $D_0^o$, the interior of
$D_0$, then $\Phi''(t)=-[(\alpha+1)p]^p\cdot p(p-1)t^{p-2}<0$, while
for $(t,y+at)\in D_1^o$ we have
\[
\Phi''(t)=-\frac{p^3(pt+p(y+at)-1)^{p-3}(1+a)}{(p-1)^{p-2}}(I_1+I_2),
\]
where
\begin{eqnarray*}
I_1 &=& pt\bigl[(p-2)(1+a)\bigl(p(\alpha+1)-1\bigr)+2\bigl(p(\alpha+1)-1-a\bigr)\bigr]\geq0,
\\
I_2 &=&\bigl(p(y+at)-1\bigr)(2\alpha+1-a)\geq0.
\end{eqnarray*}

The remaining two cases are a bit more complicated. If $(t,y+at)\in
D_2^o$, then
\[
\frac{\Phi''(t)}{Cp(1+a)^2}=J_1+J_2+J_3,
\]
where
\begin{eqnarray*}
J_1&=& h(t+y+at)^{p-2}h''(t+y+at)[h(t+y+at)-t],
\\
J_2&=& h(t+y+at)^{p-3}[h'(t+y+at)]^2[(p-1)h(t+y+at)-(p-2)t],
\\
J_3&=& -\frac{2}{a+1}h(t+y+at)^{p-2}h'(t+y+at).
\end{eqnarray*}
Now if we change $y$ and $t,$ keeping $s=t+y+at$ fixed, then
$J_1+J_2+J_3$ is a linear function of $t \in[s-1,h(s)]$. Therefore, to
prove it is non-positive, it suffices to verify this for $t=h(s)$ and
$t=s-1$. For $t=h(s)$, we have
\[
J_1+J_2+J_3=h(s)^{p-2}h'(s)\biggl[h'(s)-\frac{2}{a+1}\biggr]\leq0,
\]
since $0\leq h'(s)\leq1$ (see \eqref{backs}). If $t=s-1$, rewrite
\eqref{prod} in the form
\[
Cs\bigl(h(s)+1-s\bigr)h(s)^{p-2}h'(s)=1+C\bigl(h(s)+1-s\bigr)h(s)^{p-1}
\]
and differentiate both sides; as a result, we obtain
\begin{eqnarray*}
&&Cs\biggl[J_1+J_2+J_3+h(s)^{p-2}h'(s)\biggl(\frac{2}{a+1}-1\biggr)\biggr]
\\
&&\quad= Ch(s)^{p-2}\bigl[\bigl(h'(s)-1\bigr)h(s)+(p-2)\bigl(h(s)+1-s\bigr)h'(s)\bigr].
\end{eqnarray*}
As $h'\geq0$ and $2/(a+1)\geq1$, we will be done if we show the
right-hand side is non-positive. This is equivalent to
\[
h'(s)\bigl[h(s)+(p-2)\bigl(h(s)+1-s\bigr)\bigr]\leq h(s).
\]
Now use \eqref{backs} and substitute $h(s)=r$, noting that
$h(s)+1-s=1-\gamma(r)$, to obtain
\[
r+(p-2)\bigl(1-\gamma(r)\bigr)\leq r\bigl(1+\gamma'(r)\bigr),
\]
or $r\gamma'(r)\geq(p-2)(1-\gamma(r))$, which is \eqref{wazne}.

Finally, suppose that $y+at>1$. For such $t$ we have $\Phi(t)=(y+at)^p
u(t/(y+at),1)$, hence, setting $X=t/(y+t), Y=y+at$, we easily check
that $\Phi''(t)$ equals
\[
Y^{p-2}[p(p-1)a^2u(X,1)+2a(p-1)(1-aX)u_x(X,1)+(1-aX)^2u_{xx}(X,1)].
\]
First let us derive the expressions for the partial derivatives. Using
\eqref{prod}, we have
\begin{eqnarray*}
u_x(X,1)&=&\frac{p}{X+1}\bigl[1+C\bigl(h(X+1)-X\bigr)h(X+1)^{p-1}\bigr]-\frac
{Cph(X+1)^{p-1}}{p-1},
\\
u_{xx}(X,1)&=&\frac{p(p-1)}{(X+1)^2}\bigl[1+C\bigl(h(X+1)-X\bigr)h(X+1)^{p-1}\bigr]
\\
&&{}-\frac{Cph(X+1)^{p-1}}{X+1}-\frac{Cph(X+1)^{p-2}h'(X+1)}{X+1}.
\end{eqnarray*}
Now it can be checked that
\[
\Phi''(t)Y^{2-p}/p=K_1+K_2+K_3,
\]
where
\begin{eqnarray*}
K_1&=&(p-1)\biggl(\frac{a+1}{X+1}\biggr)^2\bigl[1+C\bigl(h(X+1)-X\bigr)h(X+1)^{p-1}\bigr],
\\
K_2&=&-\frac{Ch(X+1)^{p-1}}{X+1}(1+2a-a^2X),
\\
K_3&=&-\biggl(\frac{1-aX}{X+1}\biggr)^2\cdot\frac
{1+C(h(X+1)-X)h(X+1)^{p-1}}{h(X+1)-X}
\\
&\leq&-\biggl(\frac{1-aX}{X+1}\biggr)^2\cdot Ch(X+1)^{p-1}.
\end{eqnarray*}
We may write
\begin{eqnarray*}
K_2+K_3&\leq& -\frac
{Ch(X+1)^{p-1}}{(X+1)^2}[(1+2a-a^2X)(X+1)+(1-aX)^2]
\\
&=& -\frac{Ch(X+1)^{p-1}(a+1)}{(X+1)^2}[2+X(1-a)] \leq-\biggl(\frac
{a+1}{X+1}\biggr)^2 Ch(X+1)^{p-1},
\end{eqnarray*}
where, in the last passage, we used $a\leq1$. On the other hand, as
$h$ is non-decreasing, we have
\[
1= \frac{Ch(1)^{p}}{p-1}\leq\frac{C h(X+1)^{p-1}h(1)}{p-1}.
\]
Moreover, since $x \mapsto h(x+1)-x$ is non-increasing (see \eqref
{backs}), we have $h(X+1)-X\leq h(1)$. Combining these two facts, we obtain
\begin{eqnarray*}
K_1&\leq&(p-1)\biggl(\frac{a+1}{X+1}\biggr)^2[1+Ch(1)h(X+1)^{p-1}]
\\
&\leq& \biggl(\frac{a+1}{X+1}\biggr)^2C h(X+1)^{p-1}[h(1)+(p-1)h(1)]
\\
&\leq& \biggl(\frac{a+1}{X+1}\biggr)^2C h(X+1)^{p-1},
\end{eqnarray*}
as $ph(1)=(\alpha+1)^{-1}\leq1$. This implies $K_1+K_2+K_3\leq0$ and
completes the proof.
\end{pf*}

The final property we will need is the following.

\begin{lemma}
For any $x, y, z$ such that $x\geq0$, $z>0$ and $|y|\leq z,$ we have
%
\begin{equation}\label{deriva}
U_x(x,y,z)\leq-\alpha|U_y(x,y,z)|
\end{equation}
(if $x=0$, then $U_x$ is replaced by a right-sided derivative).
\end{lemma}

\begin{pf}
It suffices to show that for fixed $y, z$, $|y|\leq z$, and $a\in
[-\alpha,\alpha]$, the function $\Phi=\Phi_{y,z,a}\dvtx [0,\infty) \to
\mathbb{R}$ given by
$ \Phi(t)=U(t,y+at,z)$
is non-increasing.
Since $\alpha\leq1$, we know from the previous lemma that $\Phi$ is
concave. Hence all we need is $\Phi'(0+)\leq0$. By symmetry, we may
assume $y\geq0$. If $y\leq1/p$, then the derivative equals $0$; in
the remaining case, we have
\[
\Phi'(0+)= -\frac{p^2(py-1)^{p-1}}{(p-1)^{p-1}}(\alpha-a)\leq0.
\]
\upqed\end{pf}

\section{\texorpdfstring{The proof of \protect\eqref{mainin}}{The proof of (1.5)}}\label{sec4}

First let us observe that it suffices to show \eqref{mainin} for
strictly positive $\alpha$. This is an immediate consequence of the
fact that $\alpha$-strong subordination implies $\alpha'$-strong
subordination for $\alpha<\alpha'$.

Suppose $f$, $g$ are as in Theorem \ref{main}. We may restrict
ourselves to the case $\|f\|_p<\infty$. Hence, by Choi's inequality
\eqref{strongtype2}, we have $\|g\|_p<\infty$. It suffices to show
that for any $n=0, 1, 2, \ldots$ we have
\[
\mathbb{E}[(g_n^*)^p - (\alpha+1)^pp^p f_n^p ]\leq0.
\]
Clearly, we may assume that $\mathbb{P}(g_0>0)=1$, simply replacing
$f$, $g$ by $f+\varepsilon$, $g+\varepsilon$ if necessary (here
$\varepsilon$ is a small positive number). In particular, this implies
$f_0>0$ almost surely. In view of the majorization \eqref{maj}, we
will be done if we show that the expectation $\mathbb
{E}U(f_n,g_n,g_n^*)$ is
non-positive for any $n$. As a matter of fact, we will show more;
namely, that the process $(U(f_n,g_n,g_n^*)_{n\geq0})$ is a
supermartingale and $\mathbb{E}U(f_0,g_0,g_0^*)\leq0$.

To this end, fix $n\geq1$ and observe that $g^*_n\leq
|g_0|+|g_1|+\cdots+|g_n|$, so $g^*_n$ belongs to $L^p$. Thus, by Lemma
\ref{boundy} and H\"{o}lder's inequality, the variables
$U(f_n,g_n,g_n^*)$, $U_x(f_{n-1},g_{n-1},g_{n-1}^*)\,\mathrm{d}f_{n}$ and
$U_y(f_{n-1},g_{n-1},g_{n-1}^*)\,\mathrm{d}g_{n}$ are integrable. Moreover, by
definition of $U$ and the inequality \eqref{convo},
\begin{eqnarray*}
\mathbb{E}(U(f_n,g_n,g_n^*)|\mathcal{F}_{n-1})&=& \mathbb
{E}(U_n(f_n,g_n,g_{n-1}^*)|\mathcal{F}
_{n-1})
\\
&=&\mathbb{E}\bigl(U(f_{n-1}+\mathrm{d}f_n,g_{n-1}+\mathrm{d}g_n,g_{n-1}^*)|\mathcal
{F}_{n-1}\bigr)
\\
&\leq& \mathbb{E}
[U(f_{n-1},g_{n-1},g_{n-1}^*)+U_x(f_{n-1},g_{n-1},g_{n-1}^*)\,\mathrm{d}f_{n}
\\
&&\hspace*{8pt}{} +U_y(f_{n-1},g_{n-1},g_{n-1}^*)\,\mathrm{d}g_n|\mathcal{F}_{n-1}]
\\
&\leq& U(f_{n-1},g_{n-1},g_{n-1}^*).
\end{eqnarray*}
The latter inequality is the consequence of the following. By \eqref
{deriva} and the submartingale property of $f$,
\begin{eqnarray*}
\mathbb{E}(U_x(f_{n-1},g_{n-1},g_{n-1}^*)\,\mathrm{d}f_{n}|\mathcal{F}
_{n-1})&=&U_x(f_{n-1},g_{n-1},g_{n-1}^*)\mathbb{E}(\mathrm{d}f_{n}|\mathcal
{F}_{n-1})
\\[-2pt]
&\leq&-\alpha|U_y(f_{n-1},g_{n-1},g_{n-1}^*)|\mathbb{E}(\mathrm{d}f_n|\mathcal
{F}_{n-1})
\\[-2pt]
& \leq&-U_y(f_{n-1},g_{n-1},g_{n-1}^*)\mathbb{E}(\mathrm{d}g_n|\mathcal
{F}_{n-1})
\\[-2pt]
&=& -\mathbb{E}(U_y(f_{n-1},g_{n-1},g_{n-1}^*)\,\mathrm{d}g_n|\mathcal{F}_{n-1}),
\end{eqnarray*}
where the second inequality is due to $\alpha$-domination.

To complete the proof, it suffices to show that $\mathbb{E}
U(f_0,g_0,g_0^*)\leq0$. However,
$U(f_0,g_0,g_0^*)=U(f_0,g_0,g_0)=g_0^pU(f_0/g_0,1,1)$ almost surely and
the estimate follows from Corollary \ref{coro3}(ii).

\section{Sharpness}\label{sec5}

We start with inequality \eqref{mainin000} and restrict ourselves to
the case when $g$ is a $\pm1$ transform of $f$. Suppose the best
constant in this estimate equals $\beta>0$. This implies the existence
of a function $W\dvtx \mathbb{R}\times\mathbb{R}\times[0,\infty)\times
[0,\infty) \to
\mathbb{R},$ which satisfies the following properties:
%
\begin{eqnarray}\label{no0}
W(1,1,1,1)&\leq&0,
\\[-2pt]
\label{no1}
W(x,y,z,w)&=&W(x,y,|x|\vee z,|y| \vee w),\qquad  \mbox{if }x,y \in\mathbb{R},
w,z\geq0,
\\[-2pt]
\label{no2}
(|y|\vee w)^p-\beta^p (|x|\vee z)^p&\leq& W(x,y,z,w),\qquad  \mbox{if
}x,y \in\mathbb{R}, w,z\geq0
\end{eqnarray}
and, furthermore,
\begin{eqnarray} \label{no4}
&&a W(x+t_1,y+\varepsilon t_1,z,w)+(1-a) W(x+t_2,y+\varepsilon t_2,z,w)
\leq W(x,y,z,w)\nonumber
\\[-10pt]\\[-10pt]
&&\quad \mbox{for any $|x|\leq z, |y| \leq w$, $\varepsilon\in\{-1,1\}$, $a
\in(0,1)$ and $t_1,\ t_2$
with $a t_1+(1-a) t_2=0$.}\qquad\nonumber
\end{eqnarray}
Indeed, one puts
%
\begin{equation}\label{5.45}
W(x,y,z,w)=\sup\{\mathbb{E}(g_n^*\vee w)^p-\beta^p\mathbb
{E}(f_n^*\vee z)^p\},
\end{equation}
where the supremum is taken over all integers $n$ and all martingales
$f$, $g$ satisfying $\mathbb{P}((f_0,g_0)=(x,y))=1$ and $\mathrm{d}f_k=\pm
\,\mathrm{d}g_k$, $k=1, 2, \ldots$ (see \cite{B1} for details). This formula
allows us to assume that $W$ is homogeneous:
$W(tx,ty,tz,tw)=tW(x,y,z,w)$ for all $x,y \in\mathbb{R}$, $z,w \geq
0$ and $t>0$.

Now the idea is to exploit the above properties of $W$ to get $\beta
\geq p$. To this end, let $\delta$ be a small number belonging to
$(0,1/p)$. By \eqref{no4} applied to $x=0$, $y=w=1$, $z=\delta
/(1+2\delta)$, $\varepsilon=1$ and $t_1=\delta$, $t_2=-1/p$, we obtain
\begin{eqnarray}\label{cig1}
W\biggl(0,1,\frac{\delta}{1+2\delta},1\biggr)&\geq& \frac{p\delta}{1+p\delta
}W\biggl(-\frac{1}{p},1-\frac{1}{p},\frac{\delta}{1+2\delta},1\biggr)\nonumber
\\[-10pt]\\[-10pt]
&&{}+\frac{1}{1+p\delta}W\biggl(\delta,1+\delta,\frac{\delta}{1+2\delta
},1+\delta\biggr).\nonumber
\end{eqnarray}
Now, by \eqref{no1} and \eqref{no2},
%
\begin{equation}\label{cig2}
W\biggl(-\frac{1}{p},1-\frac{1}{p},\frac{\delta}{1+2\delta},1\biggr)=W\biggl(-\frac
{1}{p},1-\frac{1}{p},\frac{1}{p},1\biggr)\geq1-\biggl(\frac{\beta}{p}\biggr)^p.
\end{equation}
Furthermore, by \eqref{no1},
\[
W\biggl(\delta,1+\delta,\frac{\delta}{1+2\delta},1+\delta\biggr)=W(\delta
,1+\delta,\delta,1+\delta),
\]
which, by \eqref{no4} (with $x=z=\delta$, $y=w=1+\delta$,
$\varepsilon=-1$ and $t_1=-\delta$, $t_2=\frac{1}{p}+\delta(\frac
{1}{p}-1)$), can be bounded from below by
\[
\frac{p\delta}{1+\delta}W\biggl(\frac{1+\delta}{p},1-\frac{1}{p}+\delta
\biggl(2-\frac{1}{p}\biggr),\delta,1+\delta\biggr)+\frac{1+\delta-p\delta}{1+\delta
}W(0,1+2\delta,\delta,1+\delta).\
\]
Using \eqref{no2}, we get
\[
W\biggl(\frac{1+\delta}{p},1-\frac{1}{p}+\delta\biggl(2-\frac{1}{p}\biggr),\delta
,1+\delta\biggr)
\geq(1+\delta)^p\biggl[1-\biggl(\frac{\beta}{p}\biggr)^p\biggr].
\]
Furthermore, by \eqref{no1} and the homogeneity of $W$,
\[
W(0,1+2\delta,\delta,1+\delta)=W(0,1+2\delta,\delta,1+2\delta
)=(1+2\delta)^p W\biggl(0,1,\frac{\delta}{1+2\delta},1\biggr).
\]
Now plug all the above estimates into \eqref{cig1} to get
\begin{eqnarray}\label{niebieski}
&&W\biggl(0,1,\frac{\delta}{1+2\delta},1\biggr)\biggl[1-\frac{(1+\delta-p\delta
)(1+2\delta)^p}{(1+\delta)(1+p\delta)}\biggr]\nonumber
\\[-8pt]\\[-8pt]
&&\quad\geq\frac{p\delta}{1+p\delta}\biggl[1-\biggl(\frac{\beta}{p}\biggr)^p\biggr]\bigl(1+(1+\delta)^{p-1}\bigr).\nonumber
\end{eqnarray}
Now it follows from the definition \eqref{5.45} of $W$ that
\[
W\biggl(0,1,\frac{\delta}{1+2\delta},1\biggr) \leq W(0,1,0,1).
\]
Furthermore, one easily checks that the function
\[
F(s)=1-\frac{(1+s-ps)(1+2s)^p}{(1+s)(1+ps)}, \qquad s>-\frac{1}{p},
\]
satisfies $F(0)=F'(0)=0$. Hence
\[
1-\biggl(\frac{\beta}{p}\biggr)^p \leq\frac{W(0,1,0,1)\cdot F(\delta)\cdot
(1+p\delta)}{p\delta(1+(1+\delta)^{p-1})}
\]
and letting $\delta\to0$ yields $1-(\frac{\beta}{p})^p\leq0$, or
$\beta\geq p$.

The reasoning for the inequality \eqref{mainin0000} is essentially the
same: suppose the best constant in the estimate equals $\gamma>0$.
Introduce the function $V\dvtx [0,\infty)\times\mathbb{R}\times[0,\infty
)\times
[0,\infty) \to\mathbb{R}$ by
\[
V(x,y,z,w)=\sup\{\mathbb{E}(g_n^*\vee w)^p-\gamma^p\mathbb
{E}(f_n^*\vee z)^p\},
\]
where the supremum is taken over all integers $n$, all non-negative
submartingales $f$ and all integrable sequences $g$ satisfying $\mathbb
{P}((f_0,g_0)=(x,y))=1$ and, for $k=1, 2, \ldots,$
\[
|\mathrm{d}f_k|\geq|\mathrm{d}g_k|,\qquad\alpha\mathbb{E}(\mathrm{d}f_k|\mathcal
{F}_{k-1})\geq|\mathbb{E}(\mathrm{d}g_k|\mathcal{F}_{k-1})|
\]
with probability $1$. We see that $V$ is homogeneous and satisfies the
properties analogous to \eqref{no0}--\eqref{no4} (with obvious
changes: in \eqref{no1} and \eqref{no2} one must assume $x\geq0$; in
\eqref{no2} the number $\beta$ is replaced by $\gamma;$ and, in
\eqref{no4}, we impose $x, x+t_1, x+t_2\geq0$). In addition, there is
an extra property of $V$, which corresponds to the fact that we deal
with the inequality for submartingales:
%
\begin{equation}\label{no5}
V(x+d,y+\alpha d,z,w)\leq V(x,y,z,w),\qquad  \mbox{if }x\geq0 ,y \in
\mathbb{R}
, w,z\geq0, d\geq0.
\end{equation}
Now fix $\delta\in(0,1/p)$ and apply this property with $x=0$,
$y=w=1$, $z=\delta/(1+(\alpha+1)p)$, $d=\delta$ and then use \eqref
{no1} to obtain
\begin{eqnarray}\label{poker}
V\biggl(0,1,\frac{\delta}{1+(\alpha+1)\delta},1\biggr)&\geq& V\biggl(\delta,1+\alpha
\delta,\frac{\delta}{1+(\alpha+1)\delta},1 \biggr)\nonumber
\\[-8pt]\\[-8pt]
&=&V(\delta,1+\alpha\delta,\delta,1+\alpha\delta).\nonumber
\end{eqnarray}
Using \eqref{no1}, \eqref{no2} and \eqref{no4} as above, we have
\begin{eqnarray*}
V(\delta,1+\alpha\delta,\delta,1+\alpha\delta) &\geq&\frac
{\delta(\alpha+1)p}{1+\alpha\delta}(1+\alpha\delta)^p\biggl[1-\biggl(\frac
{\gamma}{(\alpha+1)p} \biggr)^p \biggr]
\\
&&{}+\frac{1+\alpha\delta-\delta(\alpha+1)p}{1+\alpha\delta
}\bigl(1+(\alpha+1)\delta\bigr)^pV\biggl(0,1,\frac{\delta}{1+(\alpha+1)\delta},1\biggr),
\end{eqnarray*}
which, combined with \eqref{poker}, gives
\begin{eqnarray*}
&&V\biggl(0,1,\frac{\delta}{1+(\alpha+1)\delta},1\biggr)\biggl[1-\frac{1+\alpha
\delta-\delta(\alpha+1)p}{1+\alpha\delta}\bigl(1+(\alpha+1)\delta\bigr)^p\biggr]\nonumber
\\
&&\quad \geq\delta(\alpha+1)p(1+\alpha\delta)^{p-1}\biggl[1-\biggl(\frac{\gamma
}{(\alpha+1)p} \biggr)^p \biggr].\nonumber
\end{eqnarray*}
Now it suffices to use
\[
V\biggl(0,1,\frac{\delta}{1+(\alpha+1)\delta},1\biggr)\leq V(0,1,0,1)
\]
and the fact that the function
\[
G(s)=1-\frac{1+\alpha s-s(\alpha+1)p}{1+\alpha s}\bigl(1+(\alpha
+1)s\bigr)^p,\qquad s>-1/\alpha,
\]
satisfies $G(0)=G'(0)=0$, to obtain
\[
1-\biggl(\frac{\gamma}{(\alpha+1)p} \biggr)^p\leq\frac{V(0,1,0,1)G(\delta
)}{\delta(\alpha+1)p(1+\alpha\delta)^{p-1}}.
\]
Letting $\delta\to0$ gives $1-(\frac{\gamma}{(\alpha+1)p} )^p\leq
0$, or $\gamma\geq(\alpha+1)p$. This completes the proof.

\section{Inequalities for stochastic integrals and It\^{o}
processes}\label{sec6}

In this section we present applications of the results above. Theorem
\ref{main} in the special case $\alpha=1$ yields an interesting
inequality for the stochastic integrals. Suppose $(\Omega,\mathcal
{F},\mathbb{P})$
is a complete probability space, filtered by a non-decreasing
right-continuous family $(\mathcal{F}_t)_{t\geq0}$ of sub-$\sigma
$-fields of
$\mathcal{F}$. In addition, let $\mathcal{F}_0$ contain all the
events of probability
$0$. Suppose $X=(X_t)_{t\geq0}$ is an adapted non-negative
right-continuous submartingale with left limits and let $Y$ be the It\^
{o} integral of $H$ with respect to $X$,
\[
Y_t=H_0 X_0 +\int_{(0,t]} H_s \,\mathrm{d}X_s,\qquad t\geq0.
\]
Here $H$ is a predictable process with values in $[-1,1]$. Denote
$\|X\|_p=\sup_{t\geq0} \|X_t\|_p$ and $X^*=\sup_{t\geq0} |X_t|.$ We
will establish the following extension of Theorem \ref{main}.

\begin{theorem}\label{stockh}
Under the above conditions, we have, for any $p\geq2$,
%
\begin{equation}\label{stoc}
\|Y^*\|_p \leq2p\|X\|_p,
\end{equation}
and the constant $2p$ is the best possible. It is already the best
possible in the weaker estimate
\[
\|Y^*\|_p \leq2p\|X^*\|_p.
\]
\end{theorem}

\begin{pf}
The constant $2p$ is optimal even in the discrete-time setting, so all
we need is to show \eqref{stoc}. This is a consequence of the
approximation results of Bichteler \cite{Bi}. We proceed as follows:
Consider the family $\textbf{Y}$ of all processes $Y$ of the form
%
\begin{equation}\label{stock}
Y_t=H_0X_0+\sum_{k=1}^n h_k[X_{\tau_k\wedge t}-X_{\tau_{k-1}\wedge t}],
\end{equation}
where $n$ is a positive integer, $h_k$ belongs to $[-1,1]$ and the
stopping times $\tau_k$ take only a finite number of finite values,
with $0=\tau_0\leq\tau_1\leq\cdots\leq\tau_n$. Let
\[
f=(X_{\tau_0},X_{\tau_1},\ldots,X_{\tau_n},X_{\tau_n},\ldots)
\]
and let $g$ be the transform of $f$ by $(H_0,h_1,h_2,\ldots
,h_n,0,0,\ldots)$. In virtue of Doob's optional sampling theorem, $f$
is a submartingale. Therefore, by Theorem \ref{main}, if $\tau_n\leq
t$ almost surely, then for $Y$ as in \eqref{stock},
\[
\|Y_t^*\|_p=\|g_n^*\|_p\leq2p\|f_n\|_p\leq2p\|X_t\|_p.
\]
Now we have that $X$ and $H$ satisfy the conditions of Proposition 4.1
of Bichteler \cite{Bi}. Thus by (2) of that proposition, if $Y$ is as
in the statement of the theorem above, then there is a sequence $(Y^j)$
of elements of $\textbf{Y}$ such that $\lim_{j\to\infty
}(Y^j-Y)^*=0$ almost surely. Hence, by Fatou's lemma,
\[
\|Y^*_t\|_p \leq2p\|X_t\|_p.
\]
Now take $t\to\infty$ to complete the proof.
\end{pf}

The result above can be further strengthened. Assume that $X$ is a
non-negative submartingale and $X=X_0+M+A$ stands for its Doob--Meyer
decomposition, uniquely determined by the condition that $A$ is
predictable. Let $\alpha
\in[0,1]$ be fixed and suppose $\phi$, $\psi$ are predictable
processes satisfying $|\phi_s|\leq1$ and $|\psi_s|\leq\alpha$ for
all $s$. Consider the It\^{o} process $Y$ such that $|Y_0|\leq X_0$ and
\[
Y_t=Y_0+\int_{0+}^t\phi_s\, \mathrm{d}M_s+\int_{0+}^t\psi_s \,\mathrm{d}A_s
\]
for all $t\geq0$. We have the following sharp bound.

\begin{theorem}\label{stockho}
For $X$, $Y$ as above, we have
\[
\|Y^*\|_p\leq(\alpha+1)p \|X\|_p
\]
and the inequality is sharp. So is the weaker estimate
\[
\|Y^*\|_p\leq(\alpha+1)p \|X^*\|_p.
\]
\end{theorem}

This result can be established using essentially the same approximation
arguments as above; we omit the details. We would only like to mention
here that there is an alternative way of proving Theorems \ref{stockh}
and \ref{stockho}, based on It\^{o}'s formula applied to the function
$u$ (as the function is not of class $C^2$, one needs some additional
``smoothing'' arguments to overcome this difficulty). See \cite{O3} or
\cite{S} for similar reasoning.

\section*{Acknowledgements}
This work was partially supported by MEiN Grant 1 PO3A 012 29 and the
Foundation for Polish Science.

\printhistory

\end{document}